\documentclass[11pt]{amsart}
\usepackage{amssymb}
\usepackage{amsmath}
\usepackage{amsthm}
\usepackage{latexsym}
\usepackage{amsfonts}
\usepackage{graphicx}
\usepackage{graphics}
\usepackage[T1]{pbsi}

\newtheorem{thm}{Theorem}[section]   
\newtheorem{lem}{Lemma}
\newtheorem{Def}{Definition}

\theoremstyle{definition}
\newtheorem{rem}{Remark}

\newcommand{\dis}{\displaystyle}
\newcommand{\el}{\ell}

\newcommand{\ra}{\;\rightarrow\;}

\newcommand{\ga}{\gamma }

\newcommand{\de}{\delta }

\newcommand{\De} {{\varDelta}}
\newcommand{\e}{\varepsilon }
\newcommand{\f}{\varphi}

\newcommand{\thi}{\theta }

\newcommand{\La} {{\varLambda}}

\newcommand{\la}{\lambda }
\newcommand{\mi}{\mu }

\newcommand{\ti}{\tau }

\newcommand{\oo}{\omega}

\newcommand{\R}{\mathbb{R}}
\newcommand{\Z}{\mathbb{Z}}
\newcommand{\N}{\mathbb{N}}

\newcommand{\cf}{{\mathcal{F}}}

\newcommand{\ct}{{\mathcal{T}}}

\newcommand{\cm}{{\mathcal{M}}}

\newcommand{\ld}{\ldots}

\newcommand{\sm}{\smallsetminus}

 \newcommand{\loc}{\mbox{\footnotesize loc}}
\newcommand{\hs}{\hfill$\square$}
    \newcommand{\vs}{\vspace*{0.2cm} \\}
    \newcommand{\vsp}{\vspace*{0.2cm}}

\begin{document}

\title[EXTREMAL SEQUENCES FOR THE BELLMAN FUNCTION AND APPLICATIONS]{Extremal sequences for the Bellman function of the dyadic maximal operator and applications to the Hardy operator}
\author{Eleftherios N. Nikolidakis}
\footnotetext{\hspace{-0.5cm}Keywords: Bellman function, dyadic, Hardy operator, maximal} \footnotetext{\hspace{-0.5cm}MSC Number:42B25}
\date{}
\maketitle
\noindent
{\bf Abstract:} We prove that the extremal sequences for the Bellman function of the dyadic maximal operator behave approximately as eigenfunctions of this operator for a specific eigenvalue. We use this result to prove the analogous one with respect to the Hardy operator.\bigskip\\
\noindent
\section{Introduction}\label{sec1}
\noindent

The dyadic maximal operator on $\R^n$ is a usefull tool in analysis and is defined by
\begin{eqnarray}
\hspace*{2cm}\cm_d\phi(x)=\sup\bigg\{\frac{1}{|Q|}\int_Q|\phi(y)|dy\;:x\in Q,\;Q\subseteq\R^n \;\text{in a dyadic cube}\bigg\},  \label{eq1.1}
\end{eqnarray}
for every $\phi\in L^1_{\loc}(\R^n)$, where the dyadic cubes are those formed by the grids
\[
2^{-N}\Z^n, \ \ \text{for} \ \ N=0,1,2,\ld\;.
\]
As is well known it satisfies the following weak type (1,1) inequality:
\begin{eqnarray}
|\{x\in\R^n:\cm_d\phi(x)>\;\la\}|\le\frac{1}{\la}\int_{\{\cm_d\phi>\la\}}|\phi(u)|du,  \label{eq1.2}
\end{eqnarray}
for every $\phi\in L^1(\R^n)$ and every $\la>0$.

It is easily seen that (\ref{eq1.2}) implies the following $L^p$-inequality
\begin{eqnarray}
\|\cm_d\phi\|_p\le\frac{p}{p-1}\|\phi\|_p, \ \ \label{eq1.3}
\end{eqnarray}
It is also easy to see that the weak type inequality (\ref{eq1.2}) is best possible while (\ref{eq1.3}) is also sharp.
(See \cite{1} for general martingales and \cite{15} for dyadic ones).

For the further study of the dyadic maximal operator it has been introduced the following function of two variables, defined by
\begin{eqnarray}
B_p(f,F)=\sup\bigg\{\frac{1}{|Q|}\int_Q(\cm_d\phi)^p:\;\phi\ge0,\;\frac{1}{|Q|}\int_Q\phi=f,
\frac{1}{|Q|}\int_Q\phi^p=F\bigg\},\hspace*{-2cm}  \label{eq1.4}
\end{eqnarray}
where $Q$ is a fixed dyadic cube and $0<f^p\le F$.

The function (\ref{eq1.4}), which is called the Bellman function of two variables of the dyadic maximal operator, is in fact independent of the cube $Q$ and it's value has been given in \cite{2}.\ More precisely it is proved there that
\[
B_p(f,F)=F\oo_p(f^p/F)^p,
\]
where $\oo_p:[0,1]\rightarrow\Big[1,\dfrac{p}{p-1}\Big]$ denotes the inverse function $H^{-1}_p$ of $H_p$ which is defined by
\[
H_p(z)=-(p-1)z^p+pz^{p-1}, \ \ \text{for} \ \ z\in\bigg[1,\frac{p}{p-1}\bigg].
\]
In fact this evaluation has been done in a much more general setting where the dyadic sets are now given as elements of a tree $\ct$ on a non-atomic probability space $(X,\mi)$. Then the associated dyadic maximal operator is defined by:
\begin{eqnarray}
\cm_\ct\phi(x)=\sup\bigg\{\frac{1}{\mi(I)}\int_I|\phi|d\mi:\;x\in I\in\ct\bigg\},  \label{eq1.5}
\end{eqnarray}
Additionally the inequalities (\ref{eq1.2}) and (\ref{eq1.3}) remain true and sharp in this setting. Moreover, if we define
\begin{eqnarray}
\hspace*{1.5cm}B'_{p,\ct}(f,F)=\sup\bigg\{\int_X(\cm_\ct\phi)^pd\mi:\;\phi\ge0,\;\int_X\phi d\mi=f,\int_X\phi^pd\mi=F\bigg\},  \label{eq1.6}
\end{eqnarray}
for $0<f^p\le F$, then $B'_{p,\ct}(f,F)=B_p(f,F)$. In particular the Bellman of the dyadic maximal operator is independent of the structure of the tree $\ct$.

Another approach for finding the value of $B_p(f,F)$ is given in \cite{3} where the following function of two variables has been introduced:
\begin{align}
S_p(f,F)=\sup\bigg\{\int^1_0\bigg(\frac{1}{t}\int^t_0g&\bigg)^pdt:\;g:(0,1]\rightarrow\R^+:\;\text{non-increasing},\nonumber\\
&\text{continuous and}\;\int^1_0g=f,\;\int^1_0g^p=F\bigg\}. \label{eq1.7}
\end{align}
The first step, as it can be seen in \cite{3}, is to prove that $S_p(f,F)=B_p(f,F)$.\ This can be viewed as a symmetrization principle of the dyadic maximal operator with respect to the Hardy operator.\ The second step is to prove that $S_p(f,F)$ has the expected value mentioned above.

Now the proof of the fact that $S_p=B_p$ can be given in an alternative way as can be seen in \cite{9}.\ More precisely it is proved there the following result. \vspace*{0.2cm}\\
\noindent
{\bf Theorem A.} {\em Given $g,h:(0,1]\rightarrow\R^+$ non-increasing integrable functions and a non-decreasing function $G:[0,+\infty)\rightarrow[0,+\infty)$ the following equality holds:
\begin{align*}
\sup\bigg\{\int_KG[(\cm_\ct\phi)^\ast]&h(t)dt:\;\phi\ge0, \phi^\ast=g,\;K\;\text{measurable subset of $[0,1]$ with}\\
&|K|=k\bigg\}=\int^k_0G\bigg(\frac{1}{t}\int_0^tg\bigg) h(t)dt,
\end{align*}
for any $k\in(0,1]$, where $\phi^\ast$ denotes the equimeasurable decreasing rearrangement\linebreak of $\phi$.  \hs \vspace*{0.2cm}}

It is obvious that Theorem A implies the equation $S_p=B_p$, and gives an immediate connection of the dyadic maximal operator with the Hardy operator.

An interesting question that arises now is the behaviour of the extremal sequences of functions for the quantities (\ref{eq1.6}) and (\ref{eq1.7}).\ The problem concerning (\ref{eq1.6}) has been solved in \cite{7} where it is proved the following:\vspace*{0.2cm} \\
\noindent
{\bf Theorem B.} {\em If $\phi_n:(X,\mi)\rightarrow\R^+$ be such that $\int\limits_X\phi_nd\mi=f$, $\int\limits_X\phi^p_nd\mi=F$,for every $n\in N$ then the following are equivalent

i) \;$\dis\lim_n\int\limits_X(\cm_\ct\phi_n)^pd\mi=F\oo_p(f^p/F)^p$ and

ii)\; $\dis\lim_n\int\limits_X|\cm_\ct\phi_n-c\phi_n|^pd\mi=0$, where $c=\oo_p(f^p/F)$ }.  \hs \vspace*{0.2cm}

Now it is interesting to search for the opposite problem concerning (\ref{eq1.7}).\ In fact we will prove the following:\vs
\noindent
{\bf Theorem 1.} {\em Let $g_n:(0,1]\rightarrow\R^+$ be a sequence of non-increasing functions continuous such that $\int\limits^1_0g_n(u)du=f$ and $\int\limits^1_0g_n^p(u)du=F$, for every $n\in\N$. Then the following are equivalent

i)\; $\dis\lim_n\int\limits^1_0\bigg(\dfrac{1}{t}\int\limits^t_0g_n\bigg)^pdt=F\oo_p(f^p/F)^p$

ii)\;$\dis\lim_n\int\limits^1_0\bigg|\dfrac{1}{t}\int\limits^t_0g_n-cg_n(t)\bigg|^pdt=0$ \vs
where $c=\oo_p(f^p/F)$.}  \hs \vsp

The proof is based on the proof of Theorem A and on the statement of Theorem B.

Concerning now the problem (\ref{eq1.6}) it can be easily seen that extremal functions do not exist (when the tree $\ct$ differentiates $L^1(X,\mi)$).\ That is for every $\phi\in L^p(X,\mi)$ with $\phi\ge0$ and $\int\limits_X\phi d\mi=f$, $\int\limits_X\phi^pd\mi=F$ we have the strict inequality $\int\limits_X(\cm_\ct\phi)^pd\mi<F\oo_p(f^p/F)^p$.

This is because of a self-similar property that is mentioned in \cite{8}, which states that for every extremal sequence $(\phi_n)$ for (\ref{eq1.6}) the following is true:
\begin{eqnarray}
\lim_n\frac{1}{\mi(I)}\int_I\phi_nd\mi=f \ \ \text{while} \ \ \lim_n\frac{1}{\mi(I)}\int_I\phi^p_nd\mi=F.  \label{eq1.8}
\end{eqnarray}
So, if $\phi$ is an extremal function for (\ref{eq1.6}), then we must have that $\dfrac{1}{\mi(I)}\int\limits_I\phi d\mi=f$ and $\dfrac{1}{\mi(I)}\int\limits_I\phi^pd\mi=F$ and since the tree $\ct$ differentiates $L^1(X,\mi)$ (because of (\ref{eq1.2})), then we must have that $\mi$-a.e the following equalities hold $\phi(x)=f$ and $\phi^p(x)=F$, that is $f^p=F$ which is the trivial case.

It turns out that the above doesn't hold for the extremal problem (\ref{eq1.7}).\ That is there exist extremal functions for (\ref{eq1.7}).\ We state it as: \vs
\noindent
{\bf Theorem 2.} {\em There exists unique  $g:(0,1]\rightarrow\R^+$ non-increasing and continuous with $\int\limits^1_0g(u)du=f$ and $\int\limits^1_0g^p(u)du=F$ such that
\begin{eqnarray}
\int^1_0\bigg(\frac{1}{t}\int^t_0g\bigg)^pdt=F\oo_p(f^p/F)^p.  \label{eq1.9}
\end{eqnarray}}
As it is expected due to Theorem 1, $g$ satisfies the following equality $\dfrac{1}{t}\int\limits^t_0g(u)du=\oo_p(f^p/F)g(t)$ for every $t\in(0,1]$ which gives immediately gives (\ref{eq1.9}). \hs  \vsp

After proving Theorem 2 we will be able to prove the following \vs
\noindent
{\bf Theorem 3.} {\em Let $g_n$ be as in Theorem 1. Then the following are equivalent
\begin{enumerate}
\item[i)] \;$\dis\lim_n\int\limits^1_0\bigg(\dfrac{1}{t}\int\limits^t_0g_n\bigg)^pdt=
    F\oo_p(f^p/F)^p$
\item[ii)]\; $\dis\lim_n\int\limits^1_0|g_n-g|^pdt=0$, where $g$ is the function constructed in Theorem 2. \hs
\end{enumerate}}

In this way we complete the discussion about the characterization of the extremal functions of the corresponding problem related to the Hardy operator. We
also remark that for the proof of Theorem 1 we need to fix a non-atomic probability space $(X,\mi)$ equipped with a tree structure $\ct$ which
differentiates $L^1(X,\mi)$. We use this measure space as a base in order to work there with measurable non-negative rearrangements of certain non
increasing functions on $(0,1]$.

We should also mention that the exact evaluation of (\ref{eq1.4}) for $p>1$ has been also given in \cite{10} by L. Slavin, A. Stokolos and V. Vasyunin which linked the computation of it to solving certain PDE's of the Monge-Amp\`{e}re type, and in this way they obtained an alternative proof of the results in \cite{2}. This method is different from that it is used in \cite{2} or \cite{6}. However the techniques that appear in the last two articles and the present one, give us the possibility to provide effective and powerful stability results (see for example \cite{7}).

We also remark that there are several problems in harmonic analysis were Bellman functions arise. Such problems (including the dyadic Carleson imbedding theorem and weighted inequalities) are described in \cite{10} (one can also see \cite{4} and \cite{5}) and also connections to stochastic optimal control are provided, from which it follows that the corresponding Bellman functions satisfy certain nonlinear second-order PDE's. We remark at last that the exact evaluation of a Bellman function is a difficult task and is connected with the deeper structure of the corresponding harmonic analysis problem. We mention also that until now several Bellman functions have been computed (see \cite{2}, \cite{3}, \cite{4}, \cite{5}, \cite{6}, \cite{11}, \cite{12}, \cite{13} and \cite{14}).

The paper is organized as follows:
In Section 2 we give some preliminary definitions and results.
In Section 4 we give an alternative proof of Theorem B, which is based on the proof of the evaluation of the Bellman function of two variables for
the dyadic maximal operator and which is presented in Section 3.
At last we prove Theorems 1 and 2 and 3 in Sections 5 and 6 and 7 respectively.
\section{Preliminaries}\label{sec2}
\noindent

Let $(X,\mi)$ be a non-atomic probability measure space.\ A set $\ct$ of measurable subsets of $X$ will be called a tree if it satisfies the conditions of the following
\begin{Def}\label{Def2.1}
\hspace*{-1cm}\begin{enumerate}
\item[i)] $X\in\ct$ and for every $I\in\ct$ we have that $\mi(I)>0$.
\item[ii)] For every $I\in\ct$ there corresponds a finite or countable subset $C(I)\subseteq\ct$ containing at least two elements such that
\begin{itemize}
\item[(a)] the elements of $C(I)$ are pairwise disjoint subsets of $I$
\item[(b)] $I=\cup C(I)$.
\end{itemize}
\item[iii)] $\ct=\bigcup\limits_{m\ge0}\ct_{(m)}$ where $\ct_{(0)}=\{X\}$ and $\ct_{(m+1)}=\bigcup\limits_{I\in\ct_{(m)}}C(I)$.
\item[iv)] We have that $\dis\lim_{m\ra\infty}\dis\sup_{I\in\ct_{(m)}}\mi(I)=0$.  \hs
\end{enumerate}
\end{Def}
Examples of trees are given in \cite{2}.\ The most known is the one given by the family of all dyadic subcubes of $[0,1]^n$.
The following has been proved in \cite{3}.
\begin{lem}\label{lem2.1}
For every $I\in\ct$ and every $a$ such that $0<a<1$ there exists a subfamily $\cf(I)\subseteq\ct$ consisting of pairwise disjoint subsets of $I$ such that
\[
\mi\bigg(\bigcup_{J\in\cf(I)}J\bigg)=\sum_{J\in\cf(I)}\mi(J)=(1-a)\mi(I).  \]
\hs
\end{lem}

We will also need the following fact obtained in \cite{9}.
\begin{lem}\label{lem2.2}
Let $\phi:(X,\mi)\rightarrow\R^+$ and $(A_j)_j$ a measurable partition of $X$ such that $\mi(A_j)>0$, $\forall\;j$.\ Then if $\int\limits_X\phi d\mi=f$ there exists a rearrangement of $\phi$, say $h$ $(h^\ast=\phi^\ast)$ such that $\dfrac{1}{\mi(A_j)}\int\limits_{A_j}h d\mi=f$, for every $j$. \hs
\end{lem}

Now given a tree on $(X,\mi)$ we define the associated dyadic maximal operator as follows
\[
\cm_\ct\phi(x)=\sup\bigg\{\frac{1}{\mi(I)}\int_I|\f|d\mi:\;x\in I\in\ct\bigg\},
\]
where $\phi\in L^1(X,\mi)$. We also recall the following from \cite{9}.
\begin{lem}\label{lem2.3}
Let $k\in(0,1]$ and $K$ measurable subset of $X$ with $\mi(K)=k$.\ Then the following inequality holds
\[
\int_KG[\cm_\ct\phi]d\mi\le\int^k_0G\bigg(\frac{1}{t}\int^t_0g(u)du\bigg)dt
\]
where $g\!=\!\phi^\ast$, $\phi\!\in\! L^1(X,\mi)$ and $G:[0,+\infty)\!\ra\![0,+\infty)$ is a non-decreasing function.  \hs
\end{lem}
\section{The Bellman Function of the Dyadic Maximal Operator}\label{sec3}
\noindent

In this section we provide a proof of the evaluation of the Bellman function of the dyadic maximal operators with respect to two variables $f,F$. The result appears in \cite{6} in a more general form, but we give a proof of this so as this work is complete. For this purpose we will need the following.
\setcounter{lem}{0}
\begin{lem}\label{lem3.1}
Let $\phi:(X,\mi)\rightarrow\R^+$ be such that
\[
\int_X\phi d\mi=f \quad\text{and} \quad \int_X\phi^pd\mi=F,
\]
where $0<f^p\le F$. Then
\[
\int_X(\cm_\ct\phi)^pd\mi\le F\cdot\oo_p(f^p/F)^p.
\]
\end{lem}
\begin{proof}
We consider the integral
\[
I=\int_X(\cm_\ct\phi)^pd\mi.
\]
By using Fubini's theorem we can write
\setcounter{equation}{0}
\begin{align}
I&=\int^{+\infty}_{\la=0}p\la^{p-1}\mi(\{\cm_\ct\phi>\la\})d\la \nonumber\\
&=\int^f_{\la=0}+\int^{+\infty}_{\la=f}p\la^{p-1}\mi(\{\cm_\ct\phi>\la\})d\la=
I_1+I_2,  \label{eq3.1}
\end{align}
where
\begin{align}
I_1&=\int^f_{\la=0}p\la^{p-1}\mi(\{\cm_\ct\phi>\la\})d\la \nonumber \\
&=\int^f_{\la=0}p\la^{p-1}\mi(X)d\la=\int^f_{\la=0}p\la^{p-1}d\la=f^p,  \label{eq3.2}
\end{align}
since $\cm_\ct\phi(x)\ge f$, for every $x\in X$.

$I_2$ is defined by
\[
I_2=\int^{+\infty}_{\la=f}p\la^{p-1}\mi(\{\cm_\ct\phi>\la\})d\la.
\]
By using inequality (\ref{eq1.2}) we conclude that
\begin{align*}
I_2&\le\int^{+\infty}_{\la=f}p\la^{p-1}\frac{1}{\la}\bigg(\int_{\{\cm_\ct\phi>\la\}}
\phi d\mi\bigg)d\la \\
&=\int^{+\infty}_{\la=f}p\la^{p-2}\bigg(\int_{\{\cm_\ct\phi>\la\}}\phi d\mi\bigg)d\la=\frac{p}{p-1}\int_X\phi(x)\big[\la^{p-1}\big]^{\cm_\ct\phi(x)}_{\la=f}d\mi(x),
\end{align*}
where in the last step we have used Fubini's theorem and the fact that $\cm_\ct\phi(x)\ge f$, $\forall\;x\in X$. Therefore
\begin{eqnarray}
I_2\le\frac{p}{p-1}\int_X\phi\cdot(\cm_\ct\phi)^{p-1}d\mi-\frac{p}{p-1}f^p. \label{eq3.3}
\end{eqnarray}
Thus from (\ref{eq3.1}), (\ref{eq3.2}) and (\ref{eq3.3}) we have as a consequence that
\begin{eqnarray}
I=\int_X(\cm_\ct\phi)^pd\mi\le-\frac{1}{p-1}f^p+\frac{p}{p-1}\int_X\phi\cdot(\cm_\ct\phi)^{p-1}d\mi. \label{eq3.4}
\end{eqnarray}
Using H\"{o}lder's inequality now, it is easy to see that for every $\phi$ as above the following inequality is true
\begin{eqnarray}
\int_X\phi(\cm_\ct\phi)^{p-1}d\mi\le\bigg(\int_X\phi^pd\mi\bigg)^{1/p}\cdot
\bigg(\int_X(\cm_\ct\phi)^pd\mi\bigg)^{(p-1)/p}.  \label{eq3.5}
\end{eqnarray}
By (\ref{eq3.4}) and (\ref{eq3.5}) we thus have
\begin{align}
I=\int_X(\cm_\ct\phi)^pd\mi&\le-\frac{1}{p-1}f^p+\frac{p}{p-1}\cdot F^{1/p}\cdot I^{(p-1)/p}\Rightarrow \nonumber \\
\frac{I}{F}&\le-\frac{1}{p-1}\cdot\frac{f^p}{F}+\bigg(\frac{p}{p-1}\bigg)\bigg(\frac{I}{F}\bigg)^{(p-1)/p}.
\label{eq3.6}
\end{align}
If we set now $J=\Big(\dfrac{I}{F}\Big)^{1/p}$, we have because of (\ref{eq3.6}) that
\begin{eqnarray}
J^p\le-\frac{1}{p-1}\cdot\frac{f^p}{F}+\frac{p}{p-1}J^{p-1},  \label{eq3.7}
\end{eqnarray}
We distinguish the two following cases:
\begin{enumerate}
\item[i)] $J\le1$. Then $J\le\oo_p(f^p/F)$, since $\oo_p$ takes values on $[1,p/(p-1)]$. Thus
\[
\bigg(\frac{I}{F}\bigg)^{1/p}\le\oo_p(f^p/F)\Rightarrow I\le F\oo_p(f^p/F)^p,
\]
and our result is trivial in this case.
\item[ii)] $J>1$. Then because of (\ref{eq3.7}) we conclude
\[
pJ^{p-1}-(p-1)J^p\ge\frac{f^p}{F}
\]
or that
\[
H_p(J)\ge\frac{f^p}{F}\Rightarrow J\le\oo_p\bigg(\frac{f^p}{F}\bigg),
\]
since $\oo_p=H^{-1}_p$. As a consequence we have that
\[
\int_X(\cm_\ct\phi)^pd\mi\le F\cdot\oo_p\bigg(\frac{f^p}{F}\bigg)^p,
\]
that is we derived the proof of our Lemma.
\end{enumerate}
\end{proof}

As we shall see in Section 6, for every $f,F$ fixed such that $0<f^p\le F$ and $p>1$ there exists $g:(0,1]\rightarrow\R^+$ non-increasing, continuous which satisfies $\dis\int^1_0g(u)du=f$, $\dis\int_0g^p(u)du=F$ and $\dfrac{1}{t}\dis\int^t_0g(u)du=cg(t)$, for every $t\in(0,1]$ where
\[
c=\oo_p\bigg(\frac{f^p}{F}\bigg).
\]
Thus the next Theorem is a consequence of Theorem A, and the results of this Section.
\begin{thm}\label{thm3.1}
Let $f,F$ be fixed such that $0<f^p\le F$ where $p>1$. Then the following equality is true
\begin{eqnarray}
\hspace*{0.8cm}\sup\bigg\{\int_X(\cm_\ct\phi)^pd\mi:\phi\ge0,\int_X\phi d\mi=f,\int_X
\phi^pd\mi=F\bigg\}=F\oo_p\bigg(\frac{f^p}{F}\bigg)^p.  \label{eq3.8}
\end{eqnarray}
\end{thm}
\section{Characterization of the Extremal Sequences for the Bellman Function}\label{sec4}
\noindent

In this section we will provide an alternative proof of Theorem B, different from that in \cite{7}, based on the proof of the evaluation of the Bellman function of the dyadic maximal operator, which is given in Section \ref{sec3}.\vspace*{0.2cm} \\
\noindent
{\bf Proof of Theorem B}.

i)\;$\Rightarrow$\;ii) Let $(\phi_n)_n$ be a sequence of functions $\phi_n:(X,\mi)\rightarrow\R^+$ such that $\dis\int_X\phi_nd\mi=f$, $\dis\int_X\phi^p_nd\mi=F$ for which $\dis\lim_n\dis\int_X(\cm_\ct\phi_n)^pd\mi=F\oo_p(f^p/F)^p$.

We will prove that
\[
\lim_n\int_X\mid\cm_\ct\phi_n-c\phi_n\mid^pd\mi=0,
\]
where $c=\oo_p\Big(\frac{f^p}{F}\Big)$.

By setting $\De_n=\{\cm_\ct\phi_n>c\phi_n\}$ and $\De'_n=X\setminus\De_n=\{\cm_\ct\phi_n\le c\phi_n\}$,
it is immediate to see that it is enough define
\[
I_n=\int_{\De_n}(\cm_\ct\phi_n-c\phi_n)^pd\mi \ \ \text{and} \ \ J_n=\int_{\De'_n}(c\phi_n-\cm_\ct\phi_n)^pd\mi,
\]
and then prove that $I_n$, $J_n\rightarrow0$, as $n\rightarrow\infty$.

For the evaluation of the Bellman function, as it is described in the previous section we used the following inequality:
\setcounter{equation}{0}
\begin{eqnarray}
\int_X\phi\cdot(\cm_\ct\phi)^{p-1}d\mi\le\bigg(\int_X\phi^pd\mi\bigg)^{1/p}\cdot
\bigg(\int_X(\cm_\ct\phi)^pd\mi\bigg)^{(p-1)/p},  \label{eq4.1}
\end{eqnarray}
which for our sequence $(\phi_n)_n$ must hold as an equality in the limit (we pass to a subsequence if necessary). We write this fact as
\begin{eqnarray}
\hspace*{1cm}\int_X\phi_n\cdot(\cm_\ct\phi_n)^{p-1}d\mi\approx\bigg(\int_X\phi^p_nd\mi\bigg)^{1/p}\cdot
\bigg(\int_X(\cm_\ct\phi_n)^pd\mi\bigg)^{(p-1)/p}.  \label{eq4.2}
\end{eqnarray}
Now, we are going to state and prove the following:
\setcounter{lem}{0}
\begin{lem}\label{lem4.1}
Under the above notation and hypotheses we have that:
\begin{eqnarray}
\hspace*{1cm}\int_{X_n}\phi_n(\cm_\ct\phi_n)^{p-1}d\mi\approx\bigg(\int_{X_n}\phi^p_nd\mi\bigg)^{1/p}\cdot
\bigg(\int_{X_n}(\cm_\ct\phi_n)^pd\mi\bigg)^{(p-1)/p},  \label{eq4.3}
\end{eqnarray}
where $X_n$ may be replaced either by $\De_n$ or $\De'_n$.
\end{lem}
\begin{proof}
Certainly the following inequalities hold true in view of H\"{o}lder's inequality. These are
\begin{eqnarray}
\hspace*{1.5cm}\int_{\De_n}\phi_n\cdot(\cm_\ct\phi_n)^{p-1}d\mi\le\bigg(\int_{\De_n}\phi_n^pd\mi\bigg)^{1/p}
\cdot\bigg(\int_{\De_n}(\cm_\ct\phi_n)^pd\mi\bigg)^{(p-1)/p}, \label{eq4.4}
\end{eqnarray}
and
\begin{eqnarray}
\hspace*{1.5cm}\int_{\De'_n}\phi_n(\cm_\ct\phi_n)^{p-1}d\mi\le\bigg(\int_{\De'_n}\phi^p_nd\mi
\bigg)^{1/p}\cdot\bigg(\int_{\De'_n}(\cm_\ct\phi_n)^pd\mi\bigg)^{(p-1)/p}, \label{eq4.5}
\end{eqnarray}
for any $n\in\N$. Adding them we obtain
\begin{align}
\int_X\phi_n\cdot(\cm_\ct\phi_n)^{p-1}d\mi\le&\,\bigg(\int_{\De_n}\phi^p_nd\mi\bigg)^{1/p}
\cdot\bigg(\int_{\De_n}(\cm_\ct\phi_n)^pd\mi\bigg)^{(p-1)/p} \nonumber \\
&+\,\bigg(\int_{\De'_n}\phi^p_nd\mi\bigg)^{1/p}\cdot\bigg(\int_{\De'_n}(\cm_\ct\phi_n)^p
d\mi\bigg)^{(p-1)/p},  \label{eq4.6}
\end{align}
We use now the following elementary inequality, which proof is given below.

For every $t,t'>0$, $s,s'>0$ such that
\[
t+t'=a>0 \ \ \text{and} \ \ s+s'=b>0 \ \ \text{and any} \ \ q\in(0,1),
\]
we have that
\begin{eqnarray}
t^q\cdot s^{1-q}+(t')^q\cdot(s^1)^{1-q}\le a^q\cdot b^{1-q},  \label{eq4.7}
\end{eqnarray}
Applying it for $q=1/p$ we obtain from (\ref{eq4.6}) the following inequality:
\[
\int_X\phi_n\cdot(\cm_\ct\phi_n)^{p-1}d\mi\le\bigg(\int_X\phi^p_nd\mi\bigg)^{1/p}\cdot
\bigg(\int_X(\cm_\ct\phi_n)^pd\mi\bigg)^{(p-1)/p}
\]
which in fact is an equality in the limit, because of our hypothesis. Thus, we must have equality in both (\ref{eq4.4}) and (\ref{eq4.5}) in the limit and our lemma is proved, as soon as we prove (\ref{eq4.7}).

Fix $t\in(0,a]$ and consider the function $F_t$ of the variable $s\in(0,b)$ defined by
\[
F_t(s)=t^q\cdot s^{1-q}+(a-t)^q\cdot(b-s)^{1-q}.
\]
Then
\[
F'_t(s)=(1-q)\bigg[\bigg(\frac{t}{s}\bigg)^q-\bigg(\frac{a-t}{b-s}\bigg)^q\bigg], \ \ s\in(0,b)
\]
so that $F'_t(s)>0$ for every $s\in(0,\frac{tb}{a})$, and $F'_t(s)<0$ for $s\in(\frac{tb}{a},b)$. Thus $F$ attains its maximum on the interval $[0,b]$ at the point $\frac{tb}{a}$. The result is now easily derived.
\end{proof}

We continue now with the proof of Theorem B.

Now we write
\begin{eqnarray}
\int_X(\cm_\ct\phi_n)^pd\mi=\int_{\De_n}(\cm_\ct\phi_n)^pd\mi+\int_{\De'_n}
(\cm_\ct\phi_n)^pd\mi.  \label{eq4.8}
\end{eqnarray}
We first assume that
\[
\int_{\De_n}\phi^p_nd\mi, \ \ \int_{\De'_n}\phi^p_nd\mi>0, \ \ \text{for any} \ \ n\in\N.
\]
Thus in view of H\"{o}lder's inequality, (\ref{eq4.4}), (\ref{eq4.5}) and (\ref{eq4.9}) we must have that
\begin{align}
\int_X(\cm_\ct\phi_n)^pd\mi\ge&\,\frac{\Big(\dis\int_{\De_n}\phi_n\cdot(\cm_\ct\phi_n)^{p-1}d\mi\Big)^{p/(p-1)}}
{\Big(\dis\int_{\De_n}\phi^p_nd\mi\Big)^{1/(p-1)}}\nonumber\\
&+\,
\frac{\Big(\dis\int_{\De'_n}\phi_n\cdot(\cm_\ct\phi_n)^{p-1}d\mi\Big)^{p/(p-1)}}
{\Big(\dis\int_{\De'_n}\phi_n^pd\mi\Big)^{1/(p-1)}}.  \label{eq4.9}
\end{align}
We use now H\"{o}lder's inequality in the following form:
\begin{eqnarray}
\hspace*{0.8cm}\frac{a^k}{b^{k-1}}+\frac{c^k}{d^{k-1}}\ge\frac{(a+c)^k}{(b+d)^{k-1}}, \ \ \text{for any} \ \ a,c\ge0, \ \ b,d>0, \ \ \text{where} \ \ k>1.  \label{eq4.10}
\end{eqnarray}
The above inequality is true as an equality if and only if
\[
\frac{a}{b}=\frac{c}{d}=\la, \ \ \text{for some} \ \ \la\in\R, \ \ \la\ge0.
\]

Thus in view of (\ref{eq4.10}), (\ref{eq4.9}) becomes:
\begin{eqnarray}
\int_X(\cm_\ct\phi_n)^pd\mi\ge\frac{\Big(\dis\int_X(\cm_\ct\phi_n)^{p-1}\phi_nd\mi\Big)^{(p-1)/p}}
{\Big(\dis\int_X\phi^p_nd\mi\Big)^{1/(p-1)}},  \label{eq4.11}
\end{eqnarray}
which is an equality in the limit, in view of the fact that $\phi_n$ is extremal for the Bellman function, that is $\dis\lim_n\dis\int_X(\cm_\ct\phi_n)^pd\mi=F\oo_p\Big(\dfrac{f^p}{F}\Big)^p$. From all the above we conclude, by passing if necessary to a subsequence that
\begin{eqnarray}
\hspace*{1cm}\lim_n\frac{\dis\int_{\De_n}\phi_n\cdot(\cm_\ct\phi_n)^{p-1}d\mi}{\dis\int_{\De_n}\phi^p_nd\mi}=
\lim_n\frac{\dis\int_{\De'_n}\phi_n\cdot(\cm_\ct\phi_n)^{p-1}d\mi}
{\dis\int_{\De'_n}\phi^p_nd\mi}=\la\in\R^+.  \label{eq4.12}
\end{eqnarray}
Thus, by the equality that holds in the limit in (\ref{eq4.9}), which is true because of the equality in (\ref{eq4.11}) we conclude that
\[
\la^{p/(p-1)}\lim_n\bigg[\int_{\De_n}\phi^p_nd\mi+\int_{\De'_n}\phi^p_nd\mi\bigg]=\lim_n
\int_X(\cm_\ct\phi_n)^pd\mi
\]
or that
\[
\la^{p/(p-1)}\cdot F=F\oo_p\bigg(\frac{f^p}{F}\bigg)^p\Rightarrow\la=\oo_p\bigg(\frac{f^p}{F}\bigg)^{p-1}.
\]
Thus by (\ref{eq4.12}) we conclude
\[
\int_{\De_n}\phi_n\cdot(\cm_\ct\phi_n)^{p-1}d\mi\approx\oo_p\bigg(\frac{f^p}{F}\bigg)^{p-1}\cdot
\bigg(\int_{\De_n}\phi^p_nd\mi\bigg), \ \ \text{and}
\]
\[
\int_{\De'_n}\phi_n(\cm_\ct\phi_n)^{p-1}d\mi\approx\oo_p\bigg(\frac{f^p}{F}\bigg)^{p-1}\cdot
\bigg(\int_{\De'_n}\phi^p_nd\mi\bigg).
\]
Then, because of Lemma \ref{lem4.1} we obtain that
\begin{eqnarray}
\int_{\De_n}(\cm_\ct\phi_n)^pd\mi\approx\oo_p\bigg(\frac{f^p}{F}\bigg)^p\cdot
\int_{\De_n}\phi^p_nd\mi, \ \ \text{and}  \label{eq4.13}
\end{eqnarray}
\begin{eqnarray}
\int_{\De'_n}(\cm_\ct\phi_n)^pd\mi\approx\oo_p\bigg(\frac{f^p}{F}\bigg)^p\cdot
\int_{\De'_n}\phi^p_nd\mi.  \label{eq4.14}
\end{eqnarray}
We will now need the following
\begin{lem}\label{lem4.2}
Suppose we are given $\oo_n:X_n\rightarrow\R^+$ where $X_n\subseteq X$, for $n\in\N$, and $w:X\rightarrow\R^+$ satisfying $w_n\ge w$ on $X_n$. Suppose also that
\[
\lim_n\int_{X_n}w^p_nd\mi=\lim_n\int_{X_n}w^pd\mi, \ \ \text{where} \ \ p>1.
\]
Then
\[
\lim_n\int_{X_n}(w_n-w)^pd\mi=0.
\]
\end{lem}
\begin{proof}
It is a simple matter to prove this lemma because of the following inequality.

For any $x>y>0$, $p>1$ the following holds $(x-y)^p\le x^p-y^p$. Thus
\[
\int_{X_n}(w_n-w)^pd\mi\le\int_{X_n}w^p_nd\mi-\int_{X_n}w^pd\mi\rightarrow0, \ \ \text{as} \ \ n\rightarrow\infty
\]
and the proof is complete.
\end{proof}

In view of Lemma \ref{lem4.2}, now and the definitions of $\De_n,\De'_n$, we see immediately that
\[
\int_{\De_n}(\cm_\ct\phi_n-c\phi_n)^pd\mi\rightarrow0 \ \ \text{and}
\]
\[
\int_{\De'_n}(c\phi_n-\cm_\ct\phi_n)^pd\mi\rightarrow0, \ \ \text{as} \ \ n\rightarrow\infty.
\]
As a consequence $\dis\int_X\mid\cm_\ct\phi_n-c\phi_n\mid^pd\mi\rightarrow0$, as $n\rightarrow\infty$ and our result is proved, in the case where
\begin{eqnarray}
\int_{\De_n}\phi^p_n>0 \ \ \text{and} \ \ \int_{\De'_n}\phi^p_nd\mi>0, \ \ \text{for any} \ \ n\in\N.  \label{eq4.15}
\end{eqnarray}
The same proof holds even if we have that (\ref{eq4.15}) is true for every $n\ge n_0$, for some $n_0\in\N$.

Assume now that
\[
\int_{\De'_n}\phi^p_nd\mi=0 \quad \text{for a fixed} \quad n\in\N.
\]
Since
\[
\De'_n=\{\cm_\ct\phi_n\le c\phi_n\} \ \ \text{and} \ \ \cm_\ct\phi_n(x)\ge f \ \ \text{for every} \ \ x\in X
\]
we conclude that
\[
f^p\mi(\De'_n)\le\int_{\De'_n}(\cm_\ct\phi_n)^pd\mi\le c^p\int_{\De'_n}\phi^p_n=0
\]
$\Rightarrow\mi(\De'_n)=0\Rightarrow\cm_\ct\phi_n>c\phi_n$ $\mi-a.c.$ on $X$. As a consequence, for our fixed $n\in\N$ we must have that
\[
\int_X(\cm_\ct\phi_n)^pd\mi>c^p\cdot\int_X\phi^p_nd\mi=F\cdot\oo_p(f^p/F)^p, \]
which cannot hold in view of Lemma \ref{lem3.1}.

Now suppose that for some subsequence of $(\phi_n)_n$ which we suppose without loss of generality that is the same as $(\phi_n)$, we have that %
\begin{eqnarray}
\int_{\De_n}\phi^p_nd\mi=0  \label{eq4.16}
\end{eqnarray}
Remember that $\De_n=\{\cm_\ct\phi_n>c\phi_n\}$.

Let then $x\in\{\phi_n=0\}$. Then if $x\in\De'_n$ we would have that $\cm_\ct\phi_n(x)\le c\phi_n(x)$ or that $\cm_\ct\phi_n(x)=0$, which is impossible, since $\cm_\ct\phi_n(y)\ge f$, for every $y\in X$. Thus
\[
\{\phi_n=0\}\subseteq\De_n\Rightarrow\De'_n\subseteq\{\phi_n>0\}.
\]
But from (\ref{eq4.16}) we have that $\dis\int_{\De'_n}\phi^p_nd\mi=F$, so if $\mi(\{\phi_n>0\}\setminus\De'_n)$ is positive we would obtain $\dis\int_{\{\phi_n>0\}}\phi^p_nd\mi>F$, which is impossible. Thus we have that
\[
\De'_n\subseteq\{\phi_n>0\} \ \ \text{and} \ \ \mi(\De'_n)=\mi(\{\phi_n>0\})
\]
for every $n\in\N$. Since integrals are not affected by adding or deleting a set of measure zero, we may suppose that
\begin{eqnarray}
\De'_n=\{\phi_n>0\}.  \label{eq4.17}
\end{eqnarray}
Because of Lemma \ref{lem4.1} we have that
\begin{eqnarray}
\hspace*{1.5cm}\int_{\De'_n}\phi_n(\cm_\ct\phi_n)^{p-1}\approx\bigg(\int_{\De'_n}\phi^p_nd\mi\bigg)^{1/p}
\cdot\bigg(\int_{\De'_n}(\cm_\ct\phi_n)^pd\mi\bigg)^{(p-1)/p},  \label{eq4.18}
\end{eqnarray}
Since (\ref{eq4.17}) holds we conclude by (\ref{eq4.18}) that
\begin{eqnarray}
\hspace*{1cm}\int_X\phi_n(\cm_\ct\phi_n)^{p-1}d\mi\approx F^{1/p}\bigg(\int_{\De'_n}(\cm_\ct\phi_n)^pd\mi\bigg)^{(p-1)/p},  \label{eq4.19}
\end{eqnarray}
But the next inequality is true in view of the extremality of the sequence of $(\phi_n)$ (see at the beginning of this section)
\begin{eqnarray}
\hspace*{0.5cm}\int_X\phi_n(\cm_\ct\phi_n)^{p-1}d\mi\approx F^{1/p}\cdot\bigg(\int_X(\cm_\ct\phi_n)^pd\mi\bigg)^{(p-1)/p}.  \label{eq4.20}
\end{eqnarray}
Thus
\[
\int_{\De'_n}(\cm_\ct\phi_n)^pd\mi\approx\int_X(\cm_\ct\phi_n)^pd\mi\Rightarrow
\int_{\De_n}(\cm_\ct\phi_n)^pd\mi\approx0,
\]
and since $\cm_\ct\phi_n\ge f$ on $X$ we conclude that $\mi(\De_n)\ra0$. Then
\[
\int_X\mid\cm_\ct\phi_n-c\phi_n\mid^pd\mi=\int_{\De_n}+\int_{\De'_n}\mid\cm_\ct
\phi_n-c\phi_n\mid^pd\mi=I_n+J_n.
\]
Then we proceed as follows:
$I_n=\int_{\De_n}(\cm_\ct\phi_n-c\phi_n)^pd\mi\le \int_{\De_n}(\cm_\ct\phi_n)^pd\mi-c^p\int_{\De_n}(\phi_n)^pd\mi $
in view of the elementary inequality used in the proof of Lemma 4.2. By all the above and by our hypothesis we conclude that
\[
I_n\approx0
\]
As for $J_n$, we have
\begin{align*}
J_n&=\int_{\De'_n}(c\phi_n-\cm_\ct\phi_n)^pd\mi\le c^p\int_{\De'_n}\phi^p_id\mi-\int_{\De'_n}(\cm_\ct\phi_n)^pd\mi \\
&\approx F\oo_p\bigg(\frac{f^p}{F}\bigg)^p-\int_X(\cm_\ct\phi_n)^pd\mi\approx0,
\end{align*}
since $(\phi_n)$ is extremal.

Thus, in any case we conclude Theorem B.
\section{Proof of Theorem 1}\label{sec5}
\noindent

We will prove Theorem 1 by arguing as in the proof of Theorem A and by using also Theorem B.

We begin with a sequence $(g_n)_n$ of non-increasing continuous functions $g_n:(0,1]\rightarrow\R^+$ such that $\int\limits^1_0g_n(u)du=f$ and $\int\limits^1_0g^p_n(u)du=F$ where $0<f^p\le F$.\ We set $c=\oo_p(f^p/F)$ and we suppose that $(g_n)_n$ is extremal for (\ref{eq1.7}), that is
\[
\lim_n\int^1_0\bigg(\frac{1}{t}\int^t_0g_n\bigg)^pdt=F\oo_p(f^p/F)^p=F\cdot c^p.
\]
Our aim is to prove that
\[
\lim_n\int^1_0\bigg|\frac{1}{t}\int^t_0g_n-cg_n(t)\bigg|^pdt=0.
\]
For this purpose it is enough to prove that
\setcounter{equation}{0}
\begin{eqnarray}
\int_{\{t:\frac{1}{t}\int\limits^t_0g_n>cg_n(t)\}}\bigg[\frac{1}{t}\int^t_0
g_n-cg_n(t)\bigg]^pdt=I_{1,n}\rightarrow0, \ \ \text{and}  \label{eq5.1}
\end{eqnarray}
\[
\int_{\{t:\frac{1}{t}\int\limits_0^tg_n<cg_n(t)\}}\bigg[cg_n(t)-\frac{1}{t}
\int^t_0g_n\bigg]^pdt=I_{2,n}\rightarrow0, \ \ \text{as} \ \ n\rightarrow\infty.
\]
We consider the first quantity in (\ref{eq5.1}) and similarly we work on the second.\ Set $A_n=\bigg\{t\in(0,1]:\dfrac{1}{t}\int\limits^t_0g_n>cg_n(t)\bigg\}$ so we need to prove that
\[
\int_{A_n}\bigg[\frac{1}{t}\int^t_0g_n-cg_n(t)\bigg]^pdt\rightarrow0, \ \ \text{as} \ \ n\rightarrow\infty.
\]
Since $(x-y)^p<x^p-y^p$, for $x>y>0$ and $p>1$ it is enough to prove that
\[
II_n=\int_{A_n}\bigg(\frac{1}{t}\int^t_0g_n\bigg)^pdt-c^p\int_{A_n}g^p_n\rightarrow0, \ \ n\rightarrow\infty.
\]
For each $A_n$, which is an open set of $(0,1]$ we consider it's connected components $I_{n,i}$, $i=1,2,\ld\;.$ So $A_n=\bigcup\limits^\infty_{i=1}I_{n,i}$, with $I_{n,i}$ open intervals in $(0,1]$ with $I_{n,i}\cap I_{n,j}=\emptyset$ for $i\neq j$.

Let $\e>0$.\ For every $n\in\N$ choose $i_n\in\N$ such that
\[
|III_n-III_{1,n}|<\e \ \ \text{and} \ \ |IV_n-IV_{1,n}|<\e
\]
where $III_n=\int\limits_{A_n}\Big(\dfrac{1}{t}\int\limits^t_0g_n\Big)^pdt$, $III_{1,n}=\int\limits_{F_n}\Big(\dfrac{1}{t}\int\limits^t_0g_n\Big)^pdt$, $IV_n=c^p\int\limits_{A_n}g^p_n$, $IV_{1,n}=c^p\int\limits_{F_n}g^p_n$, and $F_n=\bigcup\limits^{i_n}_{i=1}I_{n_1i}$.

It is clear that such choice of $i_n$ exists. Then $|II_n-II_{1,n}|<2\e$ where
\[
II_{1,n}=\int_{F_n}\bigg(\frac{1}{t}\int^t_0g_n\bigg)^pdt-c^p\int_{F_n}g^p_n.
\]
We need to find a $n_0\in\N$ such that $II_{1,n}<\e$, $\forall\;n\ge n_0$. Fix now a $g_n=:g$. We prove the following
\setcounter{lem}{0}
\begin{lem}\label{lem5.1}
There exists a family $\phi_a:(X,\mi)\ra\R^+$ of rearrangements of $g$ $(\phi^\ast_a=g$ for each $a\in(0,1)$) such that for each $\ga\in(0,1]$ there exists a family of measurable subsets of $X$, $S^{(\ga)}_a$ satisfying the following:
\[
\lim_{a\ra0^+}\int_{S^{(\ga)}_a}[\cm_\ct(\phi_a)]^pd\mi=\int^\ga_0\bigg(\frac{1}{t}
\int^t_0g\bigg)^pdt
\]
and $\dis\lim_{a\ra0^+}\mi(S^{(\ga)}_a)=\ga$.\ Moreover we have that $S^{(\ga)}_a\subseteq S^{(\ga')}_a$ for each $a$ $\ga<\ga'\le1$ and $a\in(0,1)$.  \hs
\end{lem}
\begin{proof}
We follow \cite{9}. Let $a\in(0,1)$. By using Lemma \ref{lem2.1} we choose for every $I\in\ct$ a family $\cf(I)\subseteq\ct$ of disjoint subsets of $I$ such that
\begin{eqnarray}
\sum_{J\in\cf(I)}\mi(J)=(1-a)\mi(I).  \label{eq5.2}
\end{eqnarray}
We define $S=S_a$ to be the smallest subset of $\ct$ such that $X\in S$ and for every $I\in S$, $\cf(I)\subseteq S$. We write for $I\in S$, $A_I=I\setminus\bigcup\limits_{J\in\cf(I)}J$. Then if $a_I=\mi(A_I)$ we have because of (\ref{eq5.2}) that $a_I=a\mi(I)$.\ It is also clear that
\[
S_a=\bigcup_{m\ge0}S_{a,(m)}, \ \ \text{where} \ \ S_{a,(0)}=\{X\} \ \ \text{and} \ \ S_{a,(m+1)}=\bigcup_{I\in S_{a,(m)}}\cf(I).
\]
We also define for $I\in S_a$, rank$(I)=r(I)$ to be the unique integer $m$ such that $I\in S_{a,(m)}$.

Additionally, we define for every $I\in S_a$ with $r(I)=m$
\[
\ga(I)=\ga_m=\frac{1}{a(1-a)^m}\int^{(1-a)^m}_{(1-a)^{m+1}}g(u)du.
\]
We also set for $I\in S_a$, $b_m(I)=\sum_{S\ni J\subseteq I\atop r(J)=r(I)+m}\mi(J)$. We easily then see inductively that $b_m(I)=(1-a)^m\mi(I)$.
It is also clear that for every $I\in S_a$, $I=\bigcup_{S_a\ni J\subseteq I}A_J$.

At last we define for every $m$ the measurable subset of $X$, $S_m=\bigcup\limits_{I\in S_{a,(m)}}I$.

Now, for each $m\ge0$, we choose $\ti^{(m)}_a:S_m\setminus S_{m+1}\ra\R$ such that
$$\big[\ti^{(m)}_a\big]^\ast=\Big(g\Big/\big((1-a)^{m+1},(1-a)^m\big]\Big)^\ast$$.

This is possible since $\mi(S_m\setminus S_{m+1})=\mi(S_m)-\mi(S_{m+1})=b_m(X)-b_{m+1}(X)=(1-a)^m-(1-a)^{m+1}=a(1-a)^m$. It is obvious now that $S_m\setminus S_{m+1}=\bigcup\limits_{I\in S_{a,(m)}}A_I$ and that
\[
\int_{S_m\setminus S_{m+1}}\ti^{(m)}_ad\mi=\int^{(1-a)^m}_{(1-a)^{m+1}}g(u)du\Rightarrow\frac{1}
{\mi(S_m\setminus S_{m+1})}\int_{S_m\setminus S_{m+1}}\ti_ad\mi=\ga_m.
\]
Using now Lemma \ref{lem2.2} we see that there exists a rearrangement of $\ti_a\Big/S_m\setminus S_{m+1}=\ti^{(m)}_a$ called $\phi^{(m)}_a$ for which $\dfrac{1}{a_I}\int\limits_{A_I}\phi^{(m)}_a=\ga_m$, for every $I\in S_{a,(m)}$.

Define now $\phi_a:X\ra\R^+$ by $\phi_a(x)=\phi^{(m)}_a(x)$, for $x\in S_m\setminus S_{m+1}$. Of course $\phi^\ast_a=g$.

Let now $I\in S_{a,(m)}$. Then
\begin{align}
Av_I(\phi_a)&=\frac{1}{\mi(I)}\int_I\phi_ad\mi=\frac{1}{\mi(I)}\sum_{S_a\ni J\subseteq I}\int_{A_J}\phi_ad\mi \nonumber\\
&=\frac{1}{\mi(I)}\sum_{\el\ge0}\sum_{S_a\ni J\subseteq I\atop r(J)=r(I)+\el}\int_{A_J}\phi_ad\mi \nonumber\\
&=\frac{1}{\mi(I)}\sum_{\el\ge0}\sum_{S_a\ni J\subseteq I}\ga_{m+\el}a_J \nonumber\\
&=\frac{1}{\mi(I)}\sum_{\el\ge0}\sum_{S_a\ni J\subseteq I}a\mi(J)\frac{1}{a(1-a)^{m+\el}}\int^{(1-a)^{m+\el}}_{(1-a)^{m+\el+1}}g(u)du \nonumber\\
&=\frac{1}{\mi(I)}\sum_{\el\ge0}\frac{1}{(1-a)^{m+\el}}\int^{(1-a)^{m+\el}}_{(1-a)^{m+\el+1}}
g(u)du\cdot\sum_{S_a\ni J\subseteq I\atop r(J)=m+\el}\mi(J) \nonumber\\
&=\frac{1}{\mi(I)}\sum_{\el\ge0}\frac{1}{(1-a)^{m+\el}}\int^{(1-a)^{m+\el}}_{A)^{m+\el+1}}
g(u)du\cdot b_\el(I)\nonumber \\
&=\frac{1}{(1-a)^m}\sum_{\el\ge0}\int^{(1-a)^{m+\el}}_{(1-a)^{m+\el+1}}g(u)du \nonumber\\
&=\frac{1}{(1-a)^m}\int^{(1-a)^m}_0g(u)du.  \label{eq5.3}
\end{align}

Now for $x\in S_m\setminus S_{m+1}$, there exists $I\in S_{a,(m)}$ such that $x\in I$ so
\begin{eqnarray}
\cm_\ct(\phi_a)(x)\ge Av_I(\phi_a)=\frac{1}{(1-a)^m}\int^{(1-a)^m}_0g(u)du=:\thi_m, \label{eq 5.4}
\end{eqnarray}
Since $\mi(S_m)=(1-a)^m$, for every $m\ge0$ we easily see from the above that we have
\[
[\cm_\ct(\phi_a)]^\ast(t)\ge\thi_m, \ \ \text{for every} \ \ t\in \big((1-a)^{m+1},(1-a)^m\big].
\]
For any $a,\ga\in(0,1]$ we now choose $m=m_a$ such that $(1-a)^{m+1}\le\ga<(1-a)^m$.\ So we have $\dis\lim_{a\ra0^+}(1-a)^{m_a}=\ga$.

Then using Lemma \ref{lem2.3} we have that
\begin{eqnarray}
\underset{a\ra0^+}{\lim\sup}\int_{\cup S_{a,(m_a)}}[\cm_\ct(\phi_a)]^pd\mi\le
\int^\ga_0\bigg(\frac{1}{t}\int^t_0g\bigg)^pdt<+\infty,  \label{eq5.5}
\end{eqnarray}
where $\cup S_{a,(m_a)}$ denotes the union of the elements of $S_{a,(m_a)}$.\ This is $S_{m_a}=\bigcup\limits_{I\in S_{a,(m_a)}}I$. This is true since $\mi(S_{m_a})\ra\ga$, as $a\ra0^+$.

Then
\begin{align}
\int_{S_{m_a}}(\cm_\ct\phi_a)^pd\mi&=\sum_{\el\ge m_a}\int_{S_\el\setminus S_{\el+1}}(\cm_\ct\phi_a)^pd\mi \nonumber\\
&\ge\sum_{\el\ge m_a}\bigg(\frac{1}{(1-a)^\el}\int^{(1-a)^\el}_0g(u)du\bigg)^p\mi(S_\el\setminus S_{\el+1}) \nonumber\\
&=\sum_{\el\ge m_a}\bigg(\frac{1}{(1-a)^\el}\int^{(1-a)^\el}_0g(u)du\bigg)^p\Big|\big((1-a)^{\el+1},(1-a)^\el
\big]\Big|,  \label{eq5.6}
\end{align}
Since $(1-a)^{m_a}\ra\ga$ and the right hand side of (\ref{eq5.6}) expresses a Riemann sum of the $\int\limits^{(1-a)^{m_a}}_0\Big(\dfrac{1}{t}\int\limits^t_0g\bigg)^pdt$ we conclude that
\begin{eqnarray}
\underset{\el\ra0^+}{\lim\sup}\int_{S_{m_a}}(\cm_\ct\phi_a)^pd\mi\ge\int^\ga_0
\bigg(\frac{1}{t}\int_0^tg\bigg)^pdt.  \label{eq5.7}
\end{eqnarray}
Then by (\ref{eq5.5}) we have equality on (\ref{eq5.7}).

We thus constructed the family $(\phi_a)_{a\in(0,1)}$, for which we easily see that if $0<\ga<\ga'\le1$ then $S^{(\ga)}_a\subseteq S^{(\ga')}_a$ for each $a\in(0,1)$. \hs
\end{proof}
\noindent
\begin{rem}\label{rem5.1}
It is not difficult to see by the proof of Lemma \ref{lem5.1} that for every $\el\in\N$, and $a\in(0,1)$ the following holds $h=g/(0,(1-a)^\el]$, where $h$ is defined by $h:=\Big(\phi_a\big/S_{a,(\el)}\Big)^\ast$ on $(0,(1-a)^\el]$.
\end{rem}

We now return to the proof of Theorem 1.

We remind that
\[
II_{1,n}=\int_{F_n}\bigg(\frac{1}{t}\int^t_0g_n\bigg)^pdt-c^p\int_{F_n}g^p_n=III_{1,n}-IV_{1,n}
\]
with $F_n=\bigcup\limits^{i_n}_{i=1}I_{n,i}=\bigcup\limits^{i_n}_{i=n}(a_{n,i_1}b_{n,i})$, which is a disjoint union. Thus
\[
III_{1,n}=\sum_n\bigg[\int^{b_{n,i}}_0\bigg(\frac{1}{t}\int^t_0g_n\bigg)^pdt-
\int^{a_{n,i}}_0\bigg(\frac{1}{t}\int^t_0g_n\bigg)^pdt\bigg].
\]
Now, for every $n\in\N$ we consider the corresponding to $g_n$, family $(\phi_{a,n})_{a\in(0,1)}$ and the respective subsets of $X$, $S^{(a_{n,i})}_{a,n}$, $S^{(b_{n,i})}_{a,n}$, $a\in(0,1)$, $i=1,2,\ld,n_i$ for which
\[
\mi\Big(S^{(a_{n,i})}_{a,n}\Big)\ra a_{n,i} \ \ \text{and} \ \ \mi\Big(S^{(b_{n,i})}_{a,n}\Big)\ra b_{n,i}, \ \ \text{as} \ \ a\ra0^+.
\]
We can also suppose that
\[
a_{n,i}<b_{n,i}\le a_{n,i+1}<b_{n,i+1}, \ \ i=1,2,\ld,i_n-1.
\]
Then we also have that
\[
S^{(a_{n,i})}_{a,n}\subseteq S^{(b_{n,i})}_{a,n}\subseteq S^{(a_{n,i+1})}_{a,n} \ \ \text{and of course}
\]
\begin{eqnarray}
\lim_{a\ra0^+}\int_{S^{(a_{n,i})}_{a,n}}[\cm_\ct(\phi_{a,n})]^pd\mi=
\int^{a_{n,i}}_0\bigg(\frac{1}{t}\int^t_0g_n\bigg)^pdt,  \label{eq5.8}
\end{eqnarray}
and similarly for the other endpoint $b_{n,i}$ of $I_{n,i}$. Therefore, by (\ref{eq5.8}) there exists for every $n\in\N$ an $a_{0,n}\in(0,1)$ such that $0<a<a_{0,n}\Rightarrow|III_{1,n}-V_n|<\dfrac{1}{n}$, where
\begin{align*}
V_n&=\sum^{i_n}_{i=1}\bigg[\int_{S^{(b_{n,i})}_{a,n}}(\cm_\ct\phi_{a,n})^pd\mi
-\int_{S^{(a_{n,i})}_{a,n}}(\cm_\ct\phi_{a,n})^pd\mi\bigg]\\
&=\int_{\La^{(a)}_n}(\cm_\ct\phi_{a,n})^pd\mi,\ \ \La^{(a)}_n=
\bigcup^{i_n}_{i=1}\big[S^{(b_{n,i})}_{a,n}\setminus S^{(a_{n,i})}_{a,n}\big].
\end{align*}
Additionally, we can suppose because of the relation
\[
\lim_{a\ra0^+}\int_X(\cm_\ct\phi_{a,n})^pdt=\int^1_0\bigg(\frac{1}{t}\int^t_0g_n\bigg)^pdt,
\ \ \text{for each} \ \ n\in\N
\]
and since $g_n$ is extremal for the problem (\ref{eq1.7}), that $a_{0,n}$ can be chosen such that for every $a\in(0,a_{0,n})$
\begin{eqnarray}
\bigg|\int_X(\cm_\ct\phi_{a,n})^pd\mi-F\oo_p(f^p/F)^p\bigg|<\frac{1}{n}, \ \ \text{for every} \ \ n\in\N. \label{eq5.9}
\end{eqnarray}
Choose $a'_n\in(0,a_n)$ and form the sequence
\[
\phi_{a'_n,n}=:\phi_n.
\]
Then, because of (\ref{eq5.9}) and since $\phi^\ast_n=g_n$ we have that $\phi_n$ is extremal for (\ref{eq1.6}).

Because of the Remark 5.1 we now have for every $\el\in\N$, each $n\in\N$ and $a\in(0,1)$, that
\[
\Big(\phi_{a,n}\big/S_{a,(\el)}\Big)^\ast:(0,\mi(S_\el)=(1-a)^\el]\ra\R^+
\]
is equal to $g_n\big/(0,(1-a)^\el]$. Since $\dis\lim_{a\ra0^+}\mi(\La^{(a)}_n)=|F_n|$, for every $n\in\N$ we can additionally suppose that $a_{0,n}$ satisfies the following
\[
\big|\mi(\La^{(a)}_n)-|F_n|\big|<\frac{1}{n}, \ \ \text{for every} \ \ a\in(0,a_{0,n})
\]
so if $\La_n=\La^{(a'_n)}_n$ we must have additionally, since $\phi_{a'_{n,n}}=\phi_n$, that
\begin{eqnarray}
\bigg|\int_{F_n}\bigg(\frac{1}{t}\int^t_0g_n\bigg)^pdt-\int_{\La_n}(\cm_\ct\phi_n)^p
d\mi\bigg|\le\frac{1}{n}  \label{eq5.10}
\end{eqnarray}
and that $\big|\mi(\La_n)-|F_n|\big|<\dfrac{1}{n}$, for every $n\in\N$.

It is also easy to see because of the above relations, the Remark 5.1 and the form of $\La_n$ (by passing to a subsequence if necessary), that
\begin{eqnarray}
\lim_n\int_{\La_n}\phi^p_n=\lim_n\int_{F_n}g^p_n.  \label{eq5.11}
\end{eqnarray}
We now take advantage of Theorem B.

Since $\phi_n$ is extremal for (\ref{eq1.6}) we must have that $\int\limits_X|\cm_\ct\phi_n-c\phi_n|^pd\mi\rightarrow0$, as $n\rightarrow\infty$ where $c=\oo_p(f^p/F)^p$. This implies:
\[
\int_{\La_n\cap\{\cm_\ct\phi_n\ge c\phi_n\}}(\cm_\ct\phi_n-c\phi_n)^pd\mi\rightarrow0, \ \ \text{as} \ \ n\rightarrow\infty \ \ \text{or}
\]
\[
\int_{\La'_n}(\cm_\ct\phi_n-c\phi_n)^pd\mi\rightarrow0, \ \ \text{as} \ \ n\rightarrow\infty, \ \ \text{where} \ \ \La'_n=\La_n\cap\{\cm_\ct\phi_n\ge c\phi_n\}.
\]
Since
\[
\bigg[\int_{\La'_n}(\cm_\ct\phi_n)^p\bigg]^{1/p}\le\bigg[\int_{\La'_n}
(\cm_\ct\phi_n-c\phi_n)^p\bigg]^{1/p}+\bigg[\int_{\La'_n}(c\phi_n)^p\bigg]^{1/p}
\]
we must have, because of the definition of $\La'_n$ and the above inequality that:
\[
\lim_n\int_{\La'_n}(\cm_\ct\phi_n)^p=c^p\lim_n\int_{\La'_n}\phi^p_n.
\]
In the same way we prove that:
\[
\lim_n\int_{\La_n\sm\La'_n}(\cm_\ct\phi_n)^p=c^p\lim_n\int_{\La_n\sm\La'_n}\phi^p_n, \ \ \text{so}
\]
\[
\lim_n\int_{\La_n}(\cm_\ct\phi_n)^pd\mi=c^p\lim_n\int_{\La_n}\phi^p_nd\mi.
\]
Because of (\ref{eq5.10}) and (\ref{eq5.11}) we have that
\[
\lim_n\int_{F_n}\bigg(\frac{1}{t}\int^t_0g_n\bigg)^pdt=\lim_nc^p\int_{F_n}g^p_n,
\]
and from the choice of $F_n$ we see that we must have that $II_n<2\e$, for $n\ge n_0$, for a suitable $n_0\in\N$. And this was our aim. \hs
\section{Uniqueness of extremal functions}\label{sec6}
\noindent

In this section we will prove that there exists unique $g_0:(0,1]\rightarrow\R^+$ continuous, with
\[
\int^1_0g_0(u)du=f, \ \ \int^1_0g^p_0(u)du=F \ \ \text{and}
\]
\[
\int^1_0\bigg(\frac{1}{t}\int^t_0g_0(u)du\bigg)^pdt=F\oo_p(f^p/F)^p.
\]
This is the statement of Theorem 2.\vspace*{0.2cm} \\
\noindent
{\bf Proof of Theorem 2.} By Theorem 1 it is obvious that if such a function $g_0$ exists, it must satisfies
\setcounter{equation}{0}
\begin{eqnarray}
\frac{1}{t}\int_0^tg_0(u)du=cg_0(t), \ \ \text{a.e on} \ \ (0,1], \ \ \text{where} \ \ c=\oo_p(f^p/F). \label{eq6.1}
\end{eqnarray}
Because of the continuity of $g_0$ we must have equality on (\ref{eq6.1}) in all $(0,1]$.

So, in order that $g_0$ satisfies (\ref{eq6.1}) we need to set $g_0(t)=kt^{-1+\frac{1}{c}}$, $t\in(0,1]$, and search for a constant $k$ (by solving the respective first order linear differential equation) such that
\[
\int_0^1g_0(u)du=f \ \ \text{and} \ \ \int_0g^p_0(u)du=F.
\]
The first equation becomes
\[
\int^1_0kt^{-1+\frac{1}{c}}dt=f\Leftrightarrow kc=f\Leftrightarrow k=f/c.
\]
So, we ask if $g_0$ for this $k$ satisfies the second equation. This is
\begin{align*}
&\int^1_0g^p_0(u)du=F\Leftrightarrow\frac{k^p}{\Big(-p+1+\frac{p}{c}\Big)}=F\Leftrightarrow
f^p/F=\bigg[(-p+1)+\frac{p}{c}\bigg]c^p\Leftrightarrow\\
&-(p-1)c^p+pc^{p-1}=f^p/F.
\end{align*}
But this is true because of the choice of $c=\oo_p(f^p/F)$ and $\oo_p=H^{-1}_p$ where
\[
H_p(z)=-(p-1)z^p+pz^{p-1}, \ \ \text{for} \ \ t\in\bigg[1,\frac{p}{p-1}\bigg].
\]
Because now of the form of $g_0:(0,1]\ra\R^+$ we have that
\[
\frac{1}{t}\int^t_0g_0(u)du=cg_0(t),\;\forall\;t\in(0,1]\Rightarrow\int^1_0
\bigg(\frac{1}{t}\int^t_0g_0(u)du\bigg)^pdu=F\oo_p(f^p/F)^p.
\]
So $g_0$ is the only extremal function in $(0,1]$.
\section{Uniqueness of extremal sequences}\label{sec7}
\noindent
We are now able to prove Theorem 3.

The direction ii)$\Rightarrow$i) is obvious from the conditions that $g$ satisfies.

We now proceed to ii)$\Rightarrow$i)

We suppose that we are given $g_n:(0,1]\ra\R^+$ non-increasing, continuous, such that $\int\limits^1_0g_n(u)du=f$, $\int^1_0g^p_n(u)du=F$ and
\[
\lim_n\int^1_0\bigg(\frac{1}{t}\int^t_0g_n(u)du\bigg)^pdt=F\oo_p(f^p/F)^p.
\]
Using Theorem 2 we conclude that
\[
\lim_n\int^1_0\bigg|\frac{1}{t}\int^t_0g_n-cg_n(t)\bigg|^pdt.
\]
Thus there exists a subsequence $(g_{k_n})_n$ such that if
\[
F_n(t)=\frac{1}{t}\int^t_0g_n-cg_n(t), \ \ t\in(0,1], \ \ n\in\N,
\]
then $F_{k_n}\rightarrow0$ almost everywhere (with respect to Lesbesgue measure). By a well known theorem in measure theory we have because of the finiteness of the measure space $[0,1]$ that $F_{k_n}\rightarrow0$ uniformly almost everywhere on $(0,1]$. This means that there exists a sequence of Lesbesgue measurable subsets of $(0,1]$, say $(H_n)_n$, such that $H_{n+1}\subseteq H_n$, $|H_n|\le\dfrac{1}{n}$ satisfying the following condition
\[
\bigg|\frac{1}{t}\int^t_0g_{k_n}-cg_{k_n}(t)\bigg|=|F_{k_n}(t)|\le\frac{1}{n}, \ \ \forall\;t\in(0,1]\setminus H_n.
\]
Additionally from the external regularity of the Lesbesgue measure, we can suppose that $H_n$ is a disjoint union of closed intervals on $(0,1]$. Let now $t,t'\in[a,1]\setminus H_{k_n}$, where $a$ is a fixed element of $(0,1]$.

Then the following hold $(c=\oo_p(f^p/F))$
\begin{align*}
|cg_{k_n}(t)-cg_{k_n}(t')|\le&\,\bigg|cg_{k_n}(t)-\frac{1}{t}\int^t_0g_{k_n}\bigg|+
\bigg|\frac{1}{t}\int^t_0g_{k_n}-\frac{1}{t'}\int^{t'}_0g_{k_n}\bigg| \\
&+\bigg|\frac{1}{t'}\int^{t'}_0g_{k_n}-cg_{k_n}(t')\bigg|=I+II+III.
\end{align*}
Then $I\le\dfrac{1}{k_n}$ since $t\notin H_{k_n}$. Similarly for III.

We look now at the second quantity II.

We may suppose that $t'>t$, so $t'=t+\de$ for some $\de>0$. Then
\setcounter{equation}{0}
\begin{align}
II&=\frac{1}{tt'}\bigg|t'\int^t_0g_{k_n}-t\int^{t'}_0g_{k_n}\bigg| \nonumber\\
&\le\frac{1}{at}\bigg|(t+\de)\int^t_0g_{k_n}-t\int^t_0g_{k_n}-t\int^{t'}_tg_ {k_n}\bigg| \nonumber \\
&=\frac{1}{at}\bigg|\de\int^t_0g_{k_n}-t\int^{t'}_tg_{k_n}\bigg| \nonumber\\
&\le\frac{\de}{a^2}f+\frac{1}{a}\int^{t'}_tg_{k_n},  \label{eq7.1}
\end{align}
where $f=\int\limits^1_0g_{k_n}$. Now by Holder's inequality we have that
\[
\int^{t'}_tg_{k_n}\le\bigg(\int^{t'}_tg^p_{k_n}\bigg)^{1/p}|t'-t|^{1-\frac{1}{p}}
=F\de^{1-\frac{1}{p}}.
\]
Thus $II\le\dfrac{\de f}{a}+\dfrac{1}{a}\de^{1-\frac{1}{p}}F$.

We consequently have that for a given $\e>0$ and $a\in(0,1)$ there exists $\de=\de_{a,\e}>0$ for which the following implication holds
\begin{eqnarray}
t,t'\in[a,1]\setminus H_{k_n},\ \ |t-t'|<\de\Rightarrow|g_{k_n}(t)-g_{k_n}(t')|<\e, \ \ \text{for every} \ \ n\in\N.\hspace*{-2cm}  \label{eq7.2}
\end{eqnarray}
Thus $(g_{k_n})_n$ has a property of type of equicontinuity on a certain set that depends on $a$. We consider now an enumeration of the rationals in $(0,1]$, let $\{q_1,q_2,\ld,q_{k_1},\ld\}=Q\cap(0,1]$.

For every $q\in Q\cap(0,1]$ we have that $(g_{k_n}(q))_n$ is a bounded sequence of real numbers, because $g_{k_n}$ is a sequence of non-negative, non-increasing functions on $(0,1]$ satisfying $\int\limits^1_0g_{k_n}=f$.

By a diagonal argument we produce a subsequence which we denote again by $g_{k_n}$ such that $g_{k_n}(q)\ra\la_q$, $n\ra\infty$ where $\la_q\in\R^+$, $q\in Q\cap(0,1]$.

Let $H=\bigcap\limits^\infty_{n=1}H_{k_n}$, which is a set of Lebesgue measure zero, and suppose that $x\in(a,1)\setminus H$. Then $x>a$, and there exist a $n_0\in\N$ such that $x\notin H_{k_{n_0}}$, so that $x\notin H_{k_n}$, $\forall\;n\ge n_0$.\ Additionally, choose a sequence $(p_k)_k$ of rationals on $(a,1)\setminus H_{k_{n_0}}$ such that $p_k\rightarrow x$.\ This is possible because the set $(a,1)\setminus H_{k_{n_0}}$ is an open set.\ Thus, we have that $p_k>a$ and $p_k\notin H_{k_n}$, $n\ge n_0$, $k\in\N$.

Let now $k_0\in\N:|p_k-x|<\de$, $\forall\;k\ge k_0$, where $\de$ is the one given in (\ref{eq7.2}).

We then have that $|g_{k_n}(x)-g_{k_m}(p_{k_0})|<\e$, for every $n\in\N$. Thus, for every such $x$, and every $n,m\in\N$ we have that
\begin{align*}
|g_{k_n}(x)-g_{k_m}(x)|\le&\,|g_{k_n}(x)-g_{k_n}(p_{k_0})|+|g_{k_n}(p_{k_0})
-g_{k_m}(p_{k_0})|  \\
&+|g_{k_m}(p_{k_0})-g_{k_m}(x)|<2\e+|g_{k_n}(p_{k_0})-g_{k_m}(p_{k_0})|.
\end{align*}
But $(g_{k_n}(p_{k_0}))_n$ is convergent sequence, thus Cauchy. Then $(g_{k_n}(x))_n$ is a Cauchy sequence for every $x\in(a,1)\setminus H$ for every $a\in(0,1]$.

Thus $(g_{k_n}(x))_n$ is a Cauchy sequence in all $(0,1]\setminus H$.

As a consequence there exists $g'_0:(0,1]\rightarrow\R^+$ such that
\begin{eqnarray}
g_{k_n}\rightarrow g'_0 \ \ \text{a.e. on} \ \ (0,1]  \label{eq7.3}
\end{eqnarray}

Additionally by using the relation
\[  \lim_n  \int\limits_0^1 |\frac{1}{t}\int\limits_0^t g_{k_n}-cg_{k_n}(t)|^pdt=0  \]
we may assume (by passing if necessary to a subsequence) that $F_{k_n}\to 0$ almost everywhere in $(0,1]$.
Thus
\begin{eqnarray}
  \frac{1}{t}\int\limits_0^t g_{k_n}-cg_{k_n}(t)  \to 0   \label{eq7.4}
\end{eqnarray}
for almost every $t\in (0,1]$.

Our aim is to prove that $\lim_n\int\limits_0^1|g_{k_n}-g_0'|^p=0$, and that $g_0'=g_0$, where $g_0$ is the function constructed in section 6,
thus giving us $g_{k_n}\stackrel{L^p}{\longrightarrow}g_0$.
By  \eqref{eq7.3} and \eqref{eq7.4} we immediately get that 
\begin{eqnarray}
\lim_n \frac{1}{t}\int\limits_0^t g_{k_n}=cg_0'(t)      \label{eq7.5}
\end{eqnarray}
for almost every $t\in (0,1]$.

Moreover, for every $t\in (0,1]$, we get  by using Fatou's Lemma 
\begin{eqnarray}
  \frac{1}{t}\int\limits_0^t g_0'= \frac{1}{t}\int\limits_0^t \lim_n g_{k_n} \le
   \liminf_n \frac{1}{t}\int\limits_0^t g_{k_n}= cg_0'(t)      \label{eq7.6}.
\end{eqnarray}
Inequality \eqref{eq7.6}  for $t=1$ gives $\int\limits_0^1g_0'\le f$, since $\int\limits_0^1g_{k_n}=f$ for all $n\in\N$.

Additionally
\[  \int\limits_0^1 (g_0')^p=\int\limits_0^1 (\lim g_{k_n})^p\le \liminf\limits_n \int\limits_0^1 g_{k_n}^p = F\]
using again Fatou's Lemma.

Now by Theorem 13.44 of \cite{16}, page 207,  we have that if $p>1$ and  $(f_n)_{n\in\mathbb{N}}$ is a sequence of nonnegative
 measurable functions 
in measure space  $(X,\mathcal{A},\mu)$  such that    
 $\sup_n \int_X f_n^p d \mu <+\infty$ 
 and $f_n\to f$ $\mu$- a.e.
then $f_n\to f$ weakly on $L^p$ that is
 $   \int_X f_n g d \mu \to   \int_X fg d\mu$ for every $g\in L^q$ where $\frac{1}{p}+\frac{1}{q}=1$.
We apply this theorem to our case:
$\sup _n \int\limits_0^1  g_{k_n}^pd \mu=  F <+\infty$,  and $g_{k_n}\to g_0'$  a.e.
 Thus
 \begin{eqnarray}
  g_{k_n}\to g_0'     \label{eq7.7}
\end{eqnarray}
weakly on $L^p$.
By \eqref{eq7.7} we immediately get 
\begin{eqnarray}
 \lim\limits_n  \int\limits_0^t  g_{k_n}=    \int\limits_0^t g_0',\quad \forall t\in (0,1].     \label{eq7.8}
\end{eqnarray}
By using now \eqref{eq7.4} and \eqref{eq7.8} we get
\begin{eqnarray}
    \frac{1}{t} \int\limits_0^t g_0' =c g_0'(t)    \label{eq7.9}
\end{eqnarray}
for almost every $t\in (0,1]$.

Thus $g_0'$ can be considered to be continuous on $(0,1]$ and \eqref{eq7.9} is true for every $t\in (0,1]$. Moreover
\eqref{eq7.8} for $t=1 $ gives

\begin{eqnarray}
   \int\limits_0^1 g_0'=f   \label{eq7.10}
\end{eqnarray}
  From \eqref{eq7.9}, which is true for every $t\in (0,1]$ and \eqref{eq7.10} we easily get that $g_0'$ has the following form
   
\begin{eqnarray}
   g_0'(t)=\frac{f}{c}t^{-1+\frac{1}{c}},\;\; t\in (0,1].   \label{eq7.11}
\end{eqnarray}
Since $c=\omega_p(\frac{f^p}{F})$ it is easily seen that $  \int\limits_0^1 (g_0')^p=F $.  Thus $g_0'=g_0$ where $g_0$ is the function 
constructed in Section 6.

Now we use the following result (see (13.47) in \cite{16}, page 208-209).
 
 If $f_n:(X,\mathcal{A},\mu)\to \R$, $n\in\N$ are such that $f_n\to f$ $\mu$-a.e. for some $f:(X,\mathcal{A},\mu)\to\R$
and $\|f_n\|_{L^1}\to \|f\|_{L^1}$, then $\|f_n-f\|_{L^1}\to 0$. By using the above mentioned theorem we get, since 
$g_{k_n}^p\to (g_0')^p$ a.e. on $(0,1]$ and 

 $\int\limits_0^1  g_{k_n}^p  = \int\limits_0^1    (g_0')^p= F  $, $\forall n\in\N$, that 
  $\int\limits_0^1  | g_{k_n}^p-(g_0')^p|\to 0$, as $n\to \infty$ which means that
   $\int\limits_0^1  | g_{k_n}^p-g_0^p|\to 0$, as $n\to \infty$.
   
 Using now the the elementary inequality $(x-y)^p\le x^p-y^p$ which holds for every $x\ge y\ge 0$, $p>1$ we immediately get that  
$\int\limits_0^1  | g_{k_n}-g_0|^p\to 0$ thus $g_{k_n}\stackrel{L^p}{\longrightarrow}g_0$, and the proof of Theorem 3 is completed.

Nikolidakis Eleftherios, University of Ioannina, Department of Mathematics, Panepistimioupolis, Ioannina, Greece.

E-mail address: lefteris@math.uoc.gr.

\end{document}